\documentclass[journal,twosided,web]{ieeecolor}
\usepackage{generic}
\usepackage{cite}
\usepackage{amsmath,amssymb,amsfonts}
\usepackage{algorithm,setspace}
\usepackage{algpseudocode}
\usepackage{graphicx}
\usepackage{textcomp}
\usepackage{dsfont}
\usepackage{color}
\usepackage{epstopdf}        
\usepackage{enumerate}
\usepackage{soul}
\usepackage{hyperref}
\usepackage{lscape}
\usepackage{multicol}
\usepackage{multirow}
\usepackage{calligra}
\usepackage{booktabs}
\usepackage{mathtools}
\usepackage{framed} 
\usepackage{picins}
\usepackage{empheq}
\usepackage{afterpage}
\usepackage{tikz}
\usepackage{graphics} 
\usepackage{epsfig} 
\usepackage{times} 
\usepackage{amsmath} 
\usepackage{amssymb}  
\usepackage{amsfonts}  
\usepackage{subfig}
\usepackage{epstopdf}
\usepackage{enumerate}
\usepackage{graphicx} 
\usepackage{xcolor}
\usepackage{comment}

\newtheorem{thm}{Theorem}
\newtheorem{prop}{Proposition}
\newtheorem{lemma}{Lemma}

\newtheorem{definition}{Definition}
\newtheorem{assumption}{Assumption}
\newtheorem{remark}{Remark}

\newcommand{\Rn}{\mathbb{R}^n}

\def\R{\mathbb{R}}

\makeatletter
\newcommand*{\rom}[1]{\expandafter\@slowromancap\romannumeral #1@}
\makeatother

\newcommand{\EE}{\mathcal{E}}

\newcommand{\G}{\mathcal{G}}

\newcommand{\bx}{{ x}}
\newcommand{\by}{{ y}}

\newcommand{\bh}{{h}}
\newcommand{\bbh}{{H}}
\newcommand{\bg}{{ g}}

\newcommand{\bu}{{\bf u}}
\newcommand{\bw}{{\bf w}}

\def\P{{\mathbb{P}}}
\newcommand{\N}{\mathcal{N}}

\def\L{{\mathcal{L}}}
\usepackage{amsmath}

\DeclareMathOperator*{\argmin}{arg\,min}

\newcommand*{\QEDB}{\hfill\ensuremath{\square}}%

\def\BibTeX{{\rm B\kern-.05em{\sc i\kern-.025em b}\kern-.08em
    T\kern-.1667em\lower.7ex\hbox{E}\kern-.125emX}}
\begin{document}
\title{Distributed Priority-Based Load Shedding over Time-Varying Communication Networks}
\author{Adel Aghajan,  Miguel Jimenez-Aparicio, Michael E. Ropp, Jorge I. Poveda
\thanks{A. Aghajan and J. I. Poveda are with the Department of Electrical and Computer Engineering, University of California, San Diego, La Jolla, CA. Miguel Jimenez-Aparicio and Michael E. Ropp are with Sandia National Laboratories, Albuquerque, NM. Corresponding author: J. I. Poveda. Email: \textit{poveda@ucsd.edu} \\ \indent This article has been co-authored by an employee of National Technology \& Engineering Solutions of Sandia, LLC under Contract No. DE-NA0003525 with the U.S. Department of Energy (DOE). The employee owns all right, title and interest in and to the article and is solely responsible for its contents. The United States Government retains and the publisher, by accepting the article for publication, acknowledges that the United States Government retains a non-exclusive, paid-up, irrevocable, world-wide license to publish or reproduce the published form of this article or allow others to do so, for United States Government purposes. The DOE will provide public access to these results of federally sponsored research in accordance with the DOE Public Access Plan \url{https://www.energy.gov/downloads/doe-public-access-plan}.}
\vspace{-0.4cm}}

\maketitle

\begin{abstract}
We study the problem of distributed optimal resource allocation on networks with actions defined on discrete spaces, with applications to adaptive under-frequency load-shedding in power systems. In this context, the primary objective is to identify an optimal subset of loads (i.e., resources) in the grid to be shed to maintain system stability whenever there is a sudden imbalance in the generation and loads. The selection of loads to be shed must satisfy demand requirements while also incorporating \emph{criticality functions} that account for socio-technical factors in the optimization process, enabling the algorithms to differentiate between network nodes with greater socio-technical value and those with less critical loads. Given the discrete nature of the state space in the optimization problem, which precludes the use of standard gradient-based approaches commonly employed in resource allocation problems with continuous action spaces, we propose a novel load-shedding algorithm based on distributed root-finding techniques and the novel concept of \emph{cumulative criticality function} (CCF). For the proposed approach, convergence conditions via Lyapunov-like techniques are established for a broad class of time-varying communication graphs that interconnect the system's regions. The theoretical results are validated through numerical examples on the Quebec 29-bus system, demonstrating the algorithm's effectiveness.
\end{abstract}

%
\vspace{-0.2cm}
\section{INTRODUCTION}
\label{sec:introduction}
\IEEEPARstart{D}{\lowercase{ecentralized}} resource allocation problems emerge in a variety of engineering systems such as power grids \cite{cherukuri2015distributed}, wireless communication systems \cite{bandwidth}, water distribution systems \cite{garcia2015modeling}, etc. In the context of energy systems, solving optimal resource allocation problems is crucial for the effective control and operation of power systems. These solutions are fundamental to modern energy infrastructure, ensuring a stable and continuous electricity supply to support societal and economic activities. However, this delicate balance can be disrupted by factors such as generation shortfalls, transmission constraints, and sudden spikes in demand, potentially resulting in system instability or even widespread blackouts \cite{muir2004blackout}, \cite{Li2007blackout}. In response to these challenges, the operational strategy of \emph{load shedding} deliberately and selectively reduces the electrical load in the grid by disconnecting a finite number of loads, aiming to mitigate risks and ensure the stability of the overall power system \cite{faranda2007proposal,concordia95isolated}.  An essential principle of effective load shedding is the prioritization of loads based on their \emph{criticality} to the system. Specifically, loads are classified according to their importance, with non-essential loads shed first to preserve power for critical infrastructure, such as hospitals and emergency services. This approach not only helps prevent system-wide blackouts but also minimizes disruptions, safeguards sensitive equipment, and enhances the grid's overall resilience. However, despite its importance in the operation of power systems, this method presents significant challenges, particularly in the context of real-time decision-making in highly heterogeneous multi-agent networks. 

\vspace{-0.3cm}
\subsection{Literature Review}
Several studies have explored effective schemes for scheduling load shedding, such as those presented in \cite{xie2020consensus,laghari2018smart,anderson1992adaptive,omar2010under}, and \cite{madiba2022under}. For example, the approaches studied in \cite{fitri2022priority} and \cite{laghari2018smart} determine load shedding by minimizing convex functions that model the cost or negative impact. However, these methods rely on pre-assigned weights to prioritize customers, which may become ineffective due to fluctuations in the energy market, new forecast data, or changes in the priority list or load shed limits. This necessitates re-tuning the weights, a non-trivial task in real-time operation. Other approaches, such as those studied in \cite{ospina2016distributed}, exploit game theoretic techniques and assume that a continuum of load can be shed at all times, while \cite{xu2011stable} studies a multi-agent system approach for load-shedding in microgrids. In fact, the most closely related work to the setting considered in this paper is presented in \cite{fitri2022priority}, which studies load-shedding problems in continuous action spaces, and under a uniform priority and criticality values for all loads within a given region. The results of \cite{fitri2022priority} also consider static network topologies. 

In practical settings, however, the optimal load-shedding problem is defined on finite action spaces, which prevents the application of standard gradient-based resource allocation techniques that rely on gradients, continuous action spaces, and convex optimization formulations \cite{cherukuri2015distributed}. In such scenarios, solving the optimal load-shedding problem, especially in distributed networked settings, becomes highly challenging. The network nodes (e.g., individual loads or groups of loads) must reach consensus on the total amount of load to be shed and \emph{how this load must be distributed across all nodes}, without having access to a centralized coordinating node. Moreover, to avoid unintended negative societal outcomes, the load-shedding process must take into account the priority of each load, whose quantifying value might be unknown to the other loads. Indeed,  in light of the growing integration of distributed energy resources (DERs), like rooftop solar and battery storage, distributed algorithms offer a promising alternative for the solution of load-shedding problems in the future, requiring less information sharing between agents and thus enhancing cybersecurity and reducing communication costs. Furthermore, distributed algorithms are recognized for their greater resilience to failures and their potential to outperform centralized methods through parallel processing. Additionally, they offer enhanced protection for sensitive data—a critical advantage in the era of distributed energy generation.

\vspace{-0.2cm}
\subsection{Contributions}
To develop a distributed algorithm that prioritizes critical loads during the shedding process, in this paper, we first introduce the concept of \emph{cumulative criticality functions} (CCFs), which represent the total load having criticality values below a given threshold. By leveraging the CCFs, we show that the load-shedding problem can be reformulated as identifying a criticality value at which the cumulative criticality function equals or exceeds the required load-shedding amount. The challenge then lies in finding this critical value in a \emph{distributed manner}. We show that this problem is equivalent to a root-finding problem in a class of time-varying functions, and we introduce a novel distributed algorithm able to recursively find this root, and therefore achieve distributed load shedding.  Finally, the efficacy of the proposed approach is validated via simulations conducted on a high-fidelity power system model of the Quebec 29-bus system.

In contrast, to existing works on resource allocation \cite{cherukuri2015distributed,barreiro2023decentralized,poveda2015shahshahani} and load-shedding problems \cite{fitri2022priority}, this paper explicitly considers the discrete nature of the load-shedding problem, wherein individual loads are typically switched on or off rather than continuously adjusted. Furthermore, we acknowledge the heterogeneity of loads within a region, recognizing that they can exhibit diverse criticality levels due to factors such as load type, function, consumer requirements, and socio-economic considerations. Additionally, our analysis considers time-varying network topologies. By accounting for these nuanced factors in our problem formulation, we improve the practicality and applicability of our proposed approach. Furthermore, while the primary emphasis is on addressing the discrete load-shedding problem, we also illustrate the algorithm's flexibility in handling continuous cases.

\vspace{-0.2cm}
\subsection{Organization}
The remainder of this paper is structured as follows: In Section \ref{sec:ProblemStatement} we present preliminaries and we formally define the problem under consideration.  Subsequently, in Section \ref{sec:mainResult}, we present our main load-shedding algorithms and convergence results. Section \ref{sec:Generalized} presents the analysis and proofs.  The performance and efficacy of the proposed algorithm are then evaluated in Section \ref{sec:simulation}. Finally, we conclude this work in Section \ref{Conclusion}. 
%
\section{PRELIMINARIES AND PROBLEM FORMULATION}\label{sec:ProblemStatement}
In this section,  we present some preliminaries and introduce the load-shedding problem studied in this paper.

\vspace{-0.2cm}
\subsection{Notation}
The following notation will be used throughout the paper. We let $[n]\triangleq\{1,\ldots,n\}$. We denote the space of real numbers by $\R$  and natural (positive integer) numbers  by $\mathbb{N}$. We denote the space of $n$-dimensional real-valued vectors by $\R^n$. In this paper, all vectors are assumed to be column vectors. The transpose of a vector $x\in \R^n$ is denoted by $x^T$.
For a vector $x\in \R^n$, $x_i$ represents the $i$th coordinate of $x$.  We denote the all-one vector in $\Rn$ by $e^n=[1,1,\ldots,1]^T$. We drop the superscript $n$ in $e^n$ whenever the dimension of the space is clear from the context.  A non-negative matrix $A$ is a row-stochastic or a column-stochastic matrix if $Ae=e$ or $e^TA= e^T$ holds, respectively. Moreover, a matrix $A$ is doubly-stochastic if it is row and column stochastic. Throughout this paper, we mainly consider undirected graphs. An undirected graph $\G=([n],\EE)$ (on $n$ vertices) is defined by a vertex set (identified by) $[n]$ and an edge set $\EE\subset [n]\times [n]$. In undirected graphs, if $(i,j)\in\EE$, then $(j,i)\in\EE$, too. The set of neighbors of vertex $j$ is represented by $\mathcal{N}_j$ for $j\in [n]$.  A graph $\G$ is connected if there exists a path between any two distinct nodes in the graph. For a matrix $A=[a_{ij}]_{n\times n}$, the associated  graph with parameter $\gamma>0$ is the graph $\G^\gamma(A)=([n],\EE^\gamma(A))$ with the edge set $\EE^\gamma(A)=\{(j,i)\mid i,j\in [n],a_{ij}>\gamma \}$. Later, we fix the value $0<\gamma<1$ throughout the paper and hence, unless otherwise stated, for notational convenience, instead of $\G^\gamma(A)$ and $\EE^\gamma(A)$, we use $\G(A)$ and $\EE(A)$, respectively.

\vspace{-0.2cm}
\subsection{The Load-Shedding Problem}
\label{sec:loadshedding}
Consider a system with $m$ electrical loads, characterized by the set 
\begin{equation}\label{setofloads}
\L\triangleq\{\ell_1, \ell_2, \dots, \ell_m\},
\end{equation}
where, for the sake of simplicity, the symbol $\ell$ can represent both a load itself and the amount of power it consumes. This means that even if two loads $\ell_i$ and $\ell_j$ have identical power demands, they are still considered distinct entities if $i\neq j$. The interpretation of $\ell$ as either the load or its power demand will be clear from the context.

We consider a scenario in which a power loss \(\P \in\mathbb{R}_{\geq0}\) occurs due to generator tripping or other causes, necessitating the shedding of \emph{at least} \(\P\) amount of load to compensate. Given the discrete nature of the set $\mathcal{L}$, shedding exactly \(\P\) may not be possible. Therefore, our goal is to identify the smallest set of loads to shed that meets or exceeds the required \(\P\). In other words, we are interested in computing a load shedding set $\mathcal{A}^*\subset \L$ that satisfies the following inclusion:
\begin{equation}\label{eqn:mainproblem}
\mathcal{A}^* \subset \argmin_{\mathcal{A}\subset\L:~\sum_{\ell \in \mathcal{A}} \ell \geq \P} ~~~\sum_{\ell \in \mathcal{A}} \ell.
\end{equation}
To guarantee feasibility for a given $\mathbb{P}\geq0$, we make the following assumption on the set of loads:
\begin{assumption}\label{assump:sumellgeP}
Given $\mathbb{P}\in\mathbb{R}_{\geq0}$, the set of loads \eqref{setofloads} satisfies $\sum_{\ell \in  \L} \ell \geq \P$.
\end{assumption}
\subsection{Incorporating Criticality Values}

\vspace{-0.2cm}
In socio-technical systems, such as power systems, certain loads  (e.g., hospital buildings) play a more critical role than others (e.g., commercial or residential buildings). Therefore, we aim to solve problem \eqref{eqn:mainproblem} by incorporating \emph{criticality values} into the load-shedding process. Specifically, to prioritize the shedding of non-critical loads within the power system, each load \(\ell\) is assigned a criticality value \(0 \leq C(\ell) \leq 1\). This value reflects the relative importance of each load, with lower values indicating less critical loads. We note that in the notation \(C(\ell)\), the symbol \(\ell\) represents the load as a distinct entity, not the quantity of the load. This means that even if two loads, \(\ell\) and \(\ell'\), have the same power demand, their criticality values may differ, i.e., \(C(\ell)\) may not be equal to \(C(\ell')\).

The objective is to integrate these criticality values into the load-shedding process to maintain system stability during sudden generation-consumption imbalances, while prioritizing the shedding of non-critical loads \(\ell_i\) with low \(C(\ell_i)\) values. Determining these criticality values requires collaboration among local authorities, energy utilities, and system operators. Furthermore, data analysis is essential to understand the impact of outages on different population demographics \cite{iliopoulos2021shedding}. 
Additionally, the criticality of loads may change over time as the duration of the outage increases.

To incorporate criticality functions into load-shedding schemes, we introduce \emph{priority-based sets}.

\vspace{0.1cm}
\begin{definition}[Priority-based Set]\label{def:priority}
A set of loads $\mathcal{A} \subseteq \L$ is deemed \emph{priority-based} if, for all $\ell \in \mathcal{A}$ and $\ell' \in \L\setminus \mathcal{A}$, the inequality ${C(\ell) \leq C(\ell')}$ holds. \QEDB  
\end{definition}
\vspace{0.1cm}

With Definition \ref{def:priority} at hand, we can formulate the priority-based load-shedding problem as identifying a set $\mathcal{A}\subset\mathcal{L}$ whose combined capacity meets or exceeds a threshold $\P$, while simultaneously adhering to the priority-based condition. Formally, this problem can be expressed as follows:
\begin{align}\label{eqn:mainProblem2}
\mathcal{A}^*\subset\argmin_{\substack{\mathcal{A} \text{ is a priority-based set} \\ \sum_{\ell\in \mathcal{A}} \ell\geq \P}}\sum_{\ell\in \mathcal{A}} \ell,
\end{align}
Note that \eqref{eqn:mainProblem2} ensures that if $C(\ell) > C(\ell')$, then load $\ell'$ will be shed before load $\ell$.
\begin{remark}
While priority-based load-shedding has been investigated before using pre-computed weights \cite{fitri2022priority,sarwar2020mixed}, the use of criticality functions in load-shedding schemes remains largely unexplored. Incorporating criticality functions offers flexibility in scenarios with highly heterogeneous loads that cannot be easily categorized into a few subsets (e.g., industrial, commercial, residential). Moreover, the criticality functions can be constructed using data-driven techniques that combine societal and technical factors. \QEDB 
\end{remark}

Under the existence of a centralized solver, the solution to problem \eqref{eqn:mainProblem2} is straightforward. In particular, assuming without loss of generality (otherwise re-order the loads) that $C(\ell_i) \leq C(\ell_{i+1})$ for $i \in [m]$, the optimal solution $\mathcal{A}^*$ can be represented as
\begin{align}\label{eqn:AiStar}
  \mathcal{A}^*=\{\ell_1,\ldots,\ell_i^*\} \text{ \ s.t. \ } \sum_{i=1}^{i^*} \ell_i\geq\P \text{ and } \sum_{i=1}^{i^*-1} \ell_i<\P.
\end{align}
On the other hand, when using a centralized solver is not feasible, computing $\mathcal{A}^*$ in a distributed manner is in general not easy. Such distributed computations are desirable in large-scale multi-area power systems, where each area might correspond to a different geographical region having access only to the information from neighboring regions, or in distributed computational systems designed to be resilient to potential failures or adversarial attacks on centralized nodes.

\vspace{-0.2cm}
\subsection{Distributed Priority-Based Load Shedding}
To solve a distributed load-shedding problem in multi-area power systems, we consider $n$ regions, and to each region $j\in\mathcal{V}:=\{1,2,\ldots,n\}$ we assign a set of loads $\L_j$, such that
\begin{align*}
    \bigcup_{j=1}^n \L_j=\L,  \text{\ \ and \ \ } \L_j\cap\L_j'=\emptyset, \ \ \forall j,j'\in [n],
\end{align*}
that is, $\{\mathcal{L}_j\}_{j=1}^n$ forms a partition of $\mathcal{L}$. Each region \(j\) has information about its own loads and their corresponding criticality values, but lacks knowledge of the loads in other regions. The criticality value of a load is influenced by two key factors: a) the nature of the load; and b) the region where the load is located. These features are explained below:
\begin{enumerate}[(a)]
\item The nature of the load, denoted \(C_n(\cdot)\), reflects the inherent importance of the load itself, such as whether it is residential, commercial, or industrial. Different types of loads naturally possess varying levels of criticality.
\item Each region has an overall criticality value, denoted \(C_r(\cdot)\), which reflects the region's relative importance within the power system. This regional criticality value can influence the criticality of individual loads within that region. 
%
\end{enumerate}
Based on this, for each $\ell\in \L_j$, we have that
\begin{align}\label{eqn:Cell}
    C(\ell)=F(C_n(\ell),C_r(j)),
\end{align}
where the function $F(\cdot,\cdot)$ is application-dependent, and where we allow each region $j$ to be aware of its own loads $\L_j$, as well as the value $C_n(\ell)$ and $C_r(j)$ for $\ell\in \L_j$, so that the criticality value \eqref{eqn:Cell} can be computed for each load in the region. 

To achieve a distributed solution to the load-shedding problem, regions must communicate with each other through a network. In this paper, we consider the scenario where this communication network can change over time, reflecting the dynamic nature of real-world communication systems. Let $\mathcal{G}(t)=([n],E(t))$ be the graph that represents the communication network at any given time $t$. Here, the $n$ regions are represented as vertices in this graph. Additionally, if there's a communication link between region $i$ and region $j$, then the edge $(i,j)$ is included in the edge set $E(t)$. We make the following ``persistent connectivity'' assumption on $\mathcal{G}(t)$:
\begin{assumption}[Communication Network]\label{assump:Network}
    There exists $B\in\mathbb{Z}_{>0}$ such that the  graph $\mathcal{G}_B(t)=([n],E_{B}(t))$, where
    \begin{align*}
       E_{B}(t)=\bigcup_{\tau=tB+1}^{(t+1)B}E(\tau),
    \end{align*}
    is  connected for all $t>0$.
\end{assumption}

To implement load shedding, it is essential to know the total required shedding amount, denoted as \(\P\). However, due to the limited information available to each region, they can only estimate a local load-shedding amount, denoted as \(p_j(t)\) for region \(j\). Without loss of generality,  we assume that, through inter-regional communication, the regions can rapidly reach a consensus on the total required shedding amount, that is
\begin{align}\label{eqn:DefW}
   \lim_{t\to\infty}\sum_{j=1}^n p_j(t)=\P.
\end{align}
This is a reasonable assumption given that, under Assumption \ref{assump:Network}, a plethora of standard distributed consensus or agreement algorithms can be used by the regions to satisfy condition \eqref{eqn:DefW}, see \cite{touri2011ergodicity,hendrickx2013convergence,touri2014endogenous,bolouki2016consensus}. Based on this, the central aim of this paper is to address the problem formulated in Equation \eqref{eqn:mainProblem2} in a distributed system characterized by a communication network that interconnects the regions and that satisfies Assumption \ref{assump:Network}.
\section{MAIN RESULTS}
\label{sec:mainResult}
In this section, we introduce the proposed distributed load-shedding algorithm for problems defined on discrete action spaces of the form \eqref{setofloads}. Before detailing the algorithm, we first introduce an equivalent formulation of the load-shedding problem \eqref{eqn:mainProblem2} using  cumulative criticality functions. We then show how this reformulated approach can be adapted to solve the load-shedding problem in a distributed manner. We also discuss extensions to continuous action spaces.

\vspace{-0.2cm}
\subsection{Cumulative Criticality Functions}
Our approach is based on the notion of \emph{cumulative criticality functions} (CCFs):

\vspace{0.1cm}
 \begin{definition}
 Consider a set of loads $\mathcal{L}$ of the form \eqref{setofloads}. We define the \emph{cumulative criticality function} (CCF) $f:\mathbb{R}\to\mathbb{R}$ of $\mathcal{L}$, as follows:
\begin{align*}
    f(z):=\sum_{\ell\in \L:~C(\ell)\leq z}\ell.
\end{align*}  
That is, $f(z)$ is the sum of the power demands of all loads whose criticality value is less than or equal to $z$.  \QEDB
\end{definition}

\vspace{0.1cm}
\begin{remark}
The CCF \(f\) can be interpreted as follows: if we are willing to shed loads with a criticality level up to and including \(z\), then the total amount of load that can be shed is represented by \(f(z)\). Note that the CCF is a right-continuous function. \QEDB 
\end{remark}

By leveraging CCFs, problem \eqref{eqn:mainProblem2} can be reformulated as finding a scalar $z^*\in\mathbb{R}$ such that
\begin{align}\label{eqn:Sol2}
    f(z^*)\geq \P \text{ \ \ and \ \ } f(z)<\P,  \ \ \forall z<z^*.
\end{align}
 We have the following Lemma, which is proved in Appendix.

\vspace{0.1cm}
\begin{lemma}\label{obs:zstar}
    The optimal threshold $z^*$ in \eqref{eqn:Sol2} corresponds to the criticality value of a specific load $\ell^* \in \mathcal{L}$, i.e., there exists $\ell^*\in\mathcal{L}$ such that $z^* = C(\ell^*)$. \QEDB
\end{lemma}

\vspace{0.1cm}
It is important to note that the solutions obtained in \eqref{eqn:AiStar} and those from \eqref{eqn:Sol2} might not be identical. Specifically, given a scalar $z^*$ satisfying \eqref{eqn:Sol2}, one cannot always guarantee that $f(z^*)=\sum_{\ell\in \mathcal{A}^*}\ell$, where $\mathcal{A}$ is given by \eqref{eqn:AiStar}. For example, consider the load quantities $\ell_1=1,~\ell_2=\ell_3=2,~\ell_4=3$ with criticality values $C(\ell_1)=0.2,~C(\ell_2)=C(\ell_3)=0.3,~C(\ell_4)=0.4$, respectively, and $\P=3$. Then, it is easy to see that $\mathcal{A}^*\subset\{\{1,2\},\{1,3\}\}$, and $z^*=0.3$, and hence 
\begin{align*}
    5=f(z^*)\not=\sum_{\ell\in \mathcal{A}^*}\ell=3.
\end{align*}
This happens because $C(\ell_2)=C(\ell_3)$. Therefore,  when using the CCF to compute $z^*$, if we shed $\ell_2$, we will also have to shed $\ell_3$, or vice versa.

However, if the quantities of criticality values of loads are \emph{distinct}, then \eqref{eqn:AiStar} and \eqref{eqn:Sol2} \emph{yield the same solution}. This is because, under the distinguishability assumption, $C(\cdot)$ is a one-to-one mapping between loads and criticality values.  More generally, we have
\begin{align*}
    f(z^*)-\sum_{\ell\in A^*}\ell\leq \sum_{\ell:C(\ell)=z^*}\ell-\min_{\ell:C(\ell)=z^*} \ell.
\end{align*}
While the solution provided by \eqref{eqn:AiStar} is the optimal solution in general, the use of CCFs and the solution given by \eqref{eqn:Sol2} are key for the proposed distributed solution of problem \eqref{eqn:mainProblem2}.

Due to their definition, it is easy to see that CCFs have the following structure:
\begin{align}\label{eqn:CCFfunction}
    f(z)=\sum_{\ell\in\L} \ell\bu\left(z-C(\ell)\right),
\end{align}
where $\bu(\cdot)$ is a right-continuous step function. Since CCFs are not Lipschitz functions, they cannot be directly used in most standard distributed algorithms for estimation and control. Therefore, we now introduce a Lipschitz approximation of \( f(\cdot) \), which is a crucial step in approximating \( z^* \).
\begin{definition}\label{def:LipCCF}
For the set of loads $\L$ given by \eqref{setofloads}, let $c\in \R_{>0}$ satisfy the following inequality:
\begin{align}\label{def:DefOfminc}
    c\leq \min_{\ell,\ell'\in\L:C(\ell)-C(\ell')>0}C(\ell)-C(\ell').
\end{align}
The function 
\begin{align*}
    \hat{f}(z)=\sum_{\ell\in\L}\ell\bw_c\left(z-C(\ell)\right),
\end{align*}
is said to be a \emph{surrogate of the CCF} \eqref{eqn:CCFfunction}, where 
\begin{align*}
    \bw_\tau(z)\triangleq
    \begin{cases}1,& \text{for }z\geq0 \\
    \frac{z}{\tau}+1,& \text{for }0>z\geq-\tau \\
    0,& \text{for }-\tau>z
    \end{cases}.
\end{align*}
Moreover, $\hat{f}$ is globally Lipschitz. \QEDB
\end{definition}
\begin{figure}
    \centering
    \includegraphics[scale=0.5]{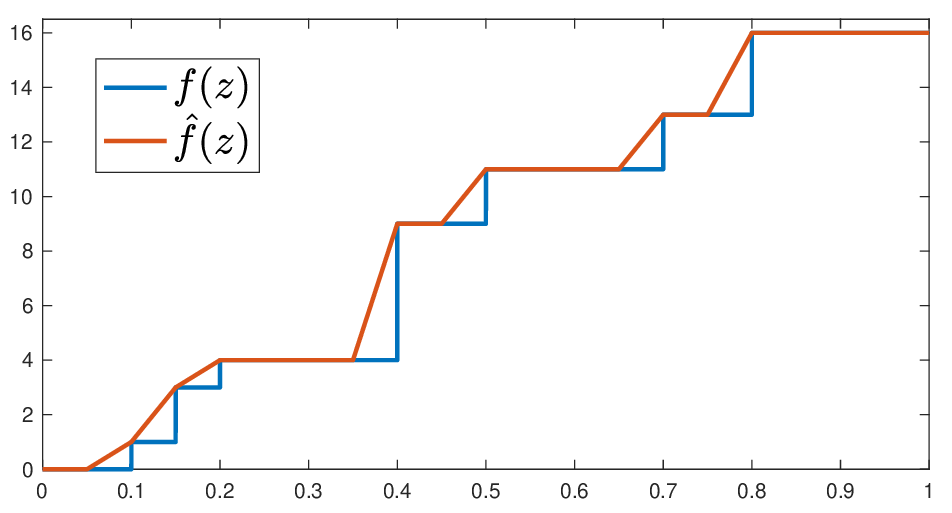}
    \caption{The  figure illustrates both the original ccf, $f(z)$, and its Lipschitz approximation, $\hat{f}(z)$, for the pair set $\{\ell\in\mathcal{L}:(\ell,C(\ell))\}$ given by $\{(1,.1), (2,.15), (1,.2), (4,.4), (1,.4), (2,.5), (2,.7), (3,.8)\}$. In this  example,  $c$ is  0.05.
}
    \label{fig:fhatf}
    \vspace{-0.5cm}
\end{figure}

It's worth noting that the parameter $c$ in \eqref{def:DefOfminc} is chosen to ensure that the original CCF \eqref{eqn:CCFfunction} and its surrogate $\hat{f}(z)$ coincide at the criticality values of all loads in the set $\mathcal{L}$. Hence, we have the following lemma, which is proved in the Appendix.

\vspace{0.1cm}
\begin{lemma}
\label{obs:fequalfhat}
    For all $\ell \in \mathcal{L}$, we have $f(C(\ell)) = \hat{f}(C(\ell))$. \QEDB
\end{lemma} 

\vspace{0.1cm}
Figure \ref{fig:fhatf} presents an example of $f(z)$ and $\hat{f}(z)$ for a pair set of loads and their corresponding criticality values $\{\ell\in\L:(\ell,C(\ell))\}$. {Given that the surrogate $\hat{f}(z)$ is a continuous function, $\hat{f}(0)=0$, and $\hat{f}(1)=\sum_{\ell\in\L}\ell\geq \P$ (Assumption \ref{assump:sumellgeP}), there exists a value $\hat{z}\in \mathbb{R}$ such that $\hat{f}(\hat{z})=\P$.} Using this value of $\hat{z}$, the following Lemma, proved in the Appendix, shows how to determine $z^*$, as defined in \eqref{eqn:Sol2}:
\begin{lemma}\label{lem:ApproxZstar}
    Let $\hat{z}\in\mathbb{R}_{\geq0}$ such that $\hat{f}(\hat{z})=\P$, and consider $z^*$ defined in \eqref{eqn:Sol2}. Then
    \begin{enumerate}
        \item if $f(z^*)>\P$, then
        $z^*=\min_{\ell: C(\ell)> \hat{z}}C(\ell),$ 
        \item if $f(z^*)=\P$, then
        $z^*=\max_{\ell: C(\ell)\leq \hat{z}}C(\ell).$
\end{enumerate}

\vspace{-0.4cm}
\QEDB
\end{lemma}

\vspace{0.1cm}
Due to the finite range of the CCF $f(\cdot)$, and the real-valued nature of the required load-shedding amount $\P$, we focus on the case $f(z^*) > \P$ in \eqref{eqn:Sol2}, since the set of points $\P$ satisfying $f(z^*) = \P$ has measure zero. 

Next, we define a local CCF and its Lipschitz surrogate for the region $j$. Specifically, for all $j\in\mathcal{V}$, we define the functions
\begin{align*}
    f_j(z)&:=\sum_{\ell\in\L_j}\ell\bu\left(z-C\left(\ell\right)\right),\\
    \hat{f}_j(z)&:=\sum_{\ell\in\L_j} \ell\bw_c\left(z-C\left(\ell\right)\right).
\end{align*}
Note that we have
\begin{align*}
    f(z)=\sum_{j=1} ^n f_j(z) \text{\ \ and \ \ } \hat{f}(z)=\sum_{j=1} ^n \hat{f}_j(z).
\end{align*}
It is assumed that each region is aware of its own CCF. Also, we assume that every region knows $c$, and hence, every region knows its own function $\hat{f}(\cdot)$. Note that regions can compute $c$ whenever they know the criticality values of all regions ($C_r(j)$ for all $j \in [n]$) and all possible load types $C_n(\ell)$. Using this information, each region can independently compute the value of $c$ via \eqref{eqn:Cell}. Note that this assumption \emph{does not} imply knowledge of the CCFs of the other regions. 

\vspace{-0.2cm}
\subsection{A Distributed Load-Shedding Algorithm}
In this subsection, we propose a distributed algorithm to compute $z^*$ in \eqref{eqn:Sol2} to solve the load-shedding problem.  The objective is for each region to determine the value of \(z^*\) and subsequently shed all loads within their region whose criticality value is less than or equal to \(z^*\). Per Lemma \ref{lem:ApproxZstar}, a preliminary task to find $z^*$ is to compute, in  a distributed way, the value of  $\hat{z}$ such that $\hat{f}(\hat{z})=\P$. The overall proposed algorithm is shown in Algorithm \ref{alg:DisLoadShed}. The following remarks are in order:
\begin{algorithm}[t!]
\setstretch{1.6}
\caption{Distributed Priority-based Load Shedding}
\begin{algorithmic}[1]
\Require $f_j(z),p_j(t)$ for all $j\in[n]$, $c$ defined in \eqref{def:DefOfminc}  
\Ensure $z^*(t)$
\State Construct $\hat{f}_j(z), \ \forall j\in[n]$ from $f_j(z)$ and $c$ 
\State Initialize $x(0)$ and $\zeta(0)$ with any numbers
\While{$\zeta(t)$ did not converge}
\For{$j\in [n]$}
    \State $x_j(t+1)\leftarrow\sum\limits_{k=1}^n w_{kj}(t)x_k(t)$
    \State \hspace{3.5cm}$+\eta(t)\left[\hat{f}_j(\bx_j(t))-p_j(t)\right]$

    \State $\zeta_j(t+1)\leftarrow\min\limits_{\ell\in \L_j :C(\ell)\geq x_j(t+1)} C(\ell)$ \label{algline:mincon1}
 
    \State Send $x_j(t+1)$ and $\zeta_j(t+1)$ to neighbors ($\N_j$) 
\EndFor
\EndWhile

\hspace{-0.7cm}Compute $z^*$ using the relationship:
\vspace{-0.4cm}
\begin{equation*}
z^* = \min_{j\in[n]} \hat{\zeta}_j.
\end{equation*}
\end{algorithmic}
\label{alg:DisLoadShed}
\end{algorithm}
\begin{enumerate}[(a)]
\item  Line 5 in Algorithm \ref{alg:DisLoadShed} can be written in vector form as follows:
\begin{align}\label{eqn:MainAlgorithm}
    {\bx}(t+1)=W(t){\bx}(t)-\eta(t)({\bg}(t)-p(t)),
\end{align}
where
\begin{itemize}
    \item $x(t)=[x_1(t),\ldots,x_n(t)]^T\in \R^n$ with $x_j(t)$ is the estimation of the region $j$ for $\hat{z}$ at time $t$,
    \item $\bx(0)\in \R$ is the initial condition,
    \item  $W(t)=[w_{kj}(t)]$ is a mixing matrix, in which if there is no communication link at time $t$ from region $j$ to region $k$, we have $w_{kj}(t)=0$,
    \item $\eta(t)$ is the step-size at time $t$, which satisfies $\lim_{t\to\infty}\eta(t)=0$. 
    \item $g(t)=[g_1(t),\ldots,g_n(t)]^T\in \R^n$ with $\bg_j(t)\triangleq \hat{f}_j(\bx_j(t))$, is a vector that contains the surrogate CCFs
    \item $p(t)=[p_1(t),\ldots,p_n(t)]^T\in \R^n$ where $p_j(t)$ is the estimation of the region $j$ for the amount of load shedding.
\end{itemize} 
\item The goal of the dynamics of the states $x_i$ is to exchange information with the other areas to converge to a common estimate of \(\hat{z}\), i.e., \(\lim_{t\to\infty} x_j(t) = \hat{z}\) for all \(j \in [n]\).
\item Assuming each region is aware of the value of $\hat{z}$, let $\hat{\zeta}_j$ represent the minimum criticality value among the loads within region $j$, i.e.,
\begin{align*}
    \hat{\zeta}_j = \min_{\ell \in \L_j:C(\ell)\geq \hat{z}} C(\ell).
\end{align*}
Note that Line 7 in Algorithm \ref{alg:DisLoadShed} can be seen as a recursive algorithm used to compute this minimum criticality value. Then, by leveraging Lemma 3, the overall minimum criticality value is given by
\begin{equation}\label{zstarequation}
z^* = \min_{j\in[n]} \hat{\zeta}_j.
\end{equation}
\item The minimum value $z^*$, given by \eqref{zstarequation}, can be computed distributively using for example a min-consensus algorithm.  Since the estimation of \(\hat{z}\) varies over time, we also can employ the ``Dynamic Min-Consensus" (DMC) protocol introduced in \cite{deplano2021dynamic}, specifically Equation (28) in \cite{deplano2021dynamic}, or its self-tuning variation \cite{deplano2022dynamicss}. 
\item Once  \(z^*\) has been computed by all areas, each area can shed all loads within their region whose criticality value is less than or equal to \(z^*\), thus solving the optimal load-shedding problem.
\end{enumerate}

To study the convergence of Algorithm 1, we make the following technical assumptions on the load-shedding estimates $p_j$ and the step size $\eta$:

\vspace{0.1cm}
\begin{assumption}\label{assump:p_jToP}
     For all $j\in[n]$, $p_j(t)$ is bounded, and $\sum_{j=1}^n p_j(t)$ converges to the total power loss $\P$. In particular, the following holds:
\begin{equation}
\left|\sum_{j=1}^n p_j(t) - \P \right| \leq \theta \eta(t)
\end{equation}
for all time steps $t$ and some positive constant $\theta$, where $\eta(t)$ is the step size in \eqref{eqn:MainAlgorithm}. \QEDB
\end{assumption}

\vspace{0.1cm}
%
%
%
%

We are now ready to present the main theoretical result of the paper:

\vspace{0.1cm}

\begin{thm}\label{thm:MainTheorem}
Assume $f(z^*)>\P$, where $z^*$ is defined in \eqref{eqn:Sol2}. Consider Algorithm \ref{alg:DisLoadShed} with inputs CCF ${f}_j$s and local load-shedding estimation $p_j(t)$s that satisfy Assumption \ref{assump:p_jToP}. Suppose that we implement Algorithm \ref{alg:DisLoadShed}  distributively on
a time-varying communication network that satisfies Assumption \ref{assump:Network}. Then, there exist a sequence 
 of matrices $\{W(t)\}$ and a sequence of step-sizes $\{\eta(t)\}$ such that
 for  $\zeta_j(t)$ defined in Algorithm~\ref{alg:DisLoadShed}
Line \ref{algline:mincon1}, we have 
\begin{align*}
   \zeta_j(t)~\text{converges  in finite time  to}  \min_{\ell\in \L_j :C(\ell)\geq \hat{z}} C(\ell),
\end{align*}
for all $j\in[n]$ and all initial conditions $\bx_j(0)\in \R$. \QEDB
\end{thm}

\vspace{0.1cm}
\begin{remark}
Drawing inspiration from distributed optimization dynamics like those in \cite{nedic2009distributed}, the proposed method addresses the challenges of priority-based load shedding in power systems without relying on gradients, convexity, or smoothness.  While existing work, such as \cite{fitri2022priority}, has explored similar problems using distributed optimization techniques, the proposed approach tackles a more realistic scenario where the set of loads is discrete and the computation of the optimal load shedding must be carried in a distributed manner. Specifically, we depart from the continuous load-shedding model with uniform priorities and static network topology employed in \cite{fitri2022priority}. Instead, we focus on discrete load shedding, where loads must be either fully shed or retained, and incorporate individual load priorities to reflect the heterogeneous nature of power systems and prioritize essential services.  Furthermore, our approach accounts for the dynamic nature of network topology, enhancing the adaptability and robustness of the load-shedding strategy. These key distinctions can enable a more nuanced and equitable response to power system disruptions. \QEDB 
\end{remark}
\vspace{0.1cm}
\begin{remark}
A common technique used in adaptive under-frequency load shedding (UFLS) relies on the swing equations of generators to compute $\P$, using the expression:
\begin{equation}\label{eqn:swingeq}
\P = - \sum_{k=1}^{N_G} \alpha_k \frac{d \omega_k(\tau)}{d \tau} \bigg{|}_{\tau=0^+},
\end{equation}
where $N_G$ is the total number of generators, $\alpha_k$ is a generator-specific coefficient based on its mechanical characteristics, and $\frac{d \omega_k(\tau)}{d \tau}$ is the derivative of the generator frequency. While recent works offer more accurate methods to estimate $\P$ via least-squares  \cite{you2021calculate}, knowledge of frequency derivatives beyond the initial time point is usually required, precluding the instantaneous estimation of the load. Thus, our goal is to demonstrate the effectiveness of our algorithm even when the load-shedding amount is not known at the beginning of the process and its continuous estimation is implemented with any mechanism that satisfies Assumption \ref{assump:p_jToP}. \QEDB 
\end{remark}

\vspace{-0.2cm}
\subsection{Variation for Continuous Action Spaces}
In the previous sections, we examined load-shedding problems where the load at each node was assumed to be confined to a finite set of values, which is common in most practical scenarios. We now show how Algorithm 1 and the concept of CCF also be easily adapted to consider settings where the load to be shed can be any arbitrary positive real number, as explored in \cite{fitri2022priority}. Specifically, the problem formulated in \cite{fitri2022priority} is as follows: There are \( n \) regions, with region \( j \) having a total load of \( L_j \) for \( j \in [n] \). The objective is to shed a real-valued amount \( l_j \), where \( 0 \leq l_j \leq L_j \), from each load \( L_j \). Each region \( j \in [n] \) has a criticality value \(C(j)=C_r(j) \), which is an integer ranging from 1 to \( n \). The goal is to solve the following problem:
\begin{align}\label{eqn:fitriproblem}
\text{Find } l_1,\ldots,&l_n \text{ s.t.} \cr
&\begin{cases}
1. \ 0\leq l_j \leq L_j, \ \forall j,\\
2. \ \sum_{j=1}^n l_j=\P \text{, and}\\ 
3. \ l_j<L_j \ \Rightarrow \ 
l_{j'}=0, \forall j: C(j')>C(j)
\end{cases} 
\end{align}
Since we are able to shed any amount of load, we have $\sum_{j=1}^n l_j=\P$, whereas in problem \eqref{eqn:mainproblem}, we do not necessarily have equality. The last condition in \eqref{eqn:fitriproblem} is similar to the priority-based set condition in Definition \ref{def:priority}.

We can extend the concept of CCF to the continuous variation of the problem. In the continuous case, the CCF is defined as:
\begin{equation}\label{eqn:CCFfunctionContinuous}
\phi(z)=\sum_{j=1}^n \phi_j(z),
\end{equation}
where $\phi_j(z)=L_j \bw_1\left(z-C(j)\right)$ for $j\in[n]$ and $\bw_\tau(z)$ is defined in \eqref{eqn:DefW}. Since $\phi(z)$ is already a Lipschitz function, there is no need for a Lipschitz surrogate.

Given that the functions \(\phi_j\) share the same structure as the functions \(f_j\), we can distributively determine the value \(\tilde{z}\) such that \(\phi(\tilde{z}) = \P\). 
 Then, a solution to \eqref{eqn:fitriproblem} can be expressed as follows:
\begin{align}\label{eqn:ContinuousSol}
\left\{\begin{array}{ll}
    l_j = L_j,& \text{if } C(j)\leq \lfloor\tilde{z}\rfloor\\
    l_j = L_j(\tilde{z}-\lfloor\tilde{z}\rfloor),& \text{if } \lfloor\tilde{z}\rfloor<C_r(j)\leq  \lceil\tilde{z}\rceil\\
    l_j =0,& \text{otherwise}
\end{array}\right.
\end{align}
The solution presented in \cite{fitri2022priority} is based on an optimization problem, which results in a solution with a small residual error. However, the CCF-based approach outlined above provides an alternative method that leverages the continuous nature of the problem and the properties of the CCF to directly determine the optimal load-shedding amounts.

\vspace{-0.2cm}
\section{ANALYSIS AND PROOFS}
\label{sec:Generalized}
In this section, we present the proofs of our results. To facilitate the analysis, we consider a broader class of functions than those explicitly defined in Equation \eqref{def:LipCCF}, which allows us to consider a generalized formulation of the load-shedding problem as a distributed iterative root-finding method.

\vspace{-0.2cm}
\subsection{Generalized Formulation and Auxiliary Results}
Consider time varying functions $h_j(z,t):\R\times\R_+\rightarrow\R$ for $j\in[n]$, a let their average limit function $H(\cdot)$ be defined as follows:
\begin{align*}
    H(z)\triangleq\frac{1}{n}\sum_{i=1}^n \lim_{t\to\infty}h_i(z,t).
\end{align*}
In the context of \eqref{eqn:MainAlgorithm}, the functions \(h_j(z,t)\) and $H(z)$ can play the role of  \(\hat{f}_j(z) - p_j(t)\) and ${(\hat{f}(z)-\P)/n}$, respectively. However, we will impose weaker assumptions on \(h_j\) and $H$, detailed below.
\begin{assumption}\label{assump:AssumpOnFuncH}
The following holds:
\begin{enumerate}
    \item The functions $h_1(z,t),\ldots,h_n(z,t)$ are uniformly bounded, i.e., there exists $M>0$ such that for all ${j\in[n],z\in\R,t\in\R}$, $|h_j(z,t)|\leq M$, and hence $H(z)\leq M$.
    \item The functions $h_j(z,t)$ for $j\in[n]$ are Lipschitz, i.e., there exists $\lambda>0$ such that for all $j\in[n]$, $z,z'\in\R$, and $t\in\R$, we have ${|h_j(z,t)-h_j(z',t)|\leq \lambda|z-z'|}$. Hence, $h_j(z,t)$ and $H(z)$ are continuous functions of $z$.
    \item There exists a $z^\star\in \R$ such that $(z-z^\star)H(z)\geq 0$ for all $z\in\R$.
\end{enumerate}

\vspace{-0.2cm}
\QEDB
\end{assumption}
Note that since $H$ is a continuous function, any $z^\star$ that satisfies item 3) in Assumption \ref{assump:AssumpOnFuncH} is a root of $H$.

 
Since our goal is to find distributed dynamics that satisfy $\lim_{t\to\infty}\bx_i(t)=z^\star$ for all $i\in[n]$, we consider a modified version of equation \eqref{eqn:MainAlgorithm},  designed for the \(h_j(z,t)\) functions. In particular, for $t\geq 0$, we consider the update rule:
\begin{align}\label{eqn:AuxAlgorithm}
    {\bx}(t+1)=W(t){\bx}(t)-\eta(t){\by}(t),
\end{align}
where, again $\bx(0)\in \R$ is the initial condition, $\bx(t)\in \R $ is the value of $x$ at time $t$, $\{\eta(t)\}$ is a step-size sequence, $\{W(t)\}$ is a mixing matrices sequence, and
\begin{align*}
    \by_j(t)\triangleq h_j(\bx_j(t),t)\in\R,~~\forall~j\in[n].
\end{align*} 
The following assumption parallels Assumption \ref{assump:p_jToP}, but is specifically tailored for the functions \(h_j(z,t)\).
\begin{assumption}\label{assump:AssumpOnConvergenceFunction}
There exist $\theta>0$ such that 
\begin{align*}
    \left|H(z)-\frac{1}{n}\sum_{i=1}^n \lim_{t\to\infty}h_i(z,t)\right|\leq \theta \eta(t).
\end{align*}
for all $i\in[n], z\in \R$. \QEDB
\end{assumption}
Since our  goal is to guarantee that  \(\lim_{t\to\infty} x_j(t) = z^\star\) for all \(j \in [n]\), sufficient information exchange between regions is needed.  The next assumption imposes some regularity conditions on the sequence of mixing matrices \(\{W(t)\}\):
\begin{assumption}\label{assump:AssumpOnConnectivity}
The following holds for all $t \geq 0$:
\begin{enumerate}[(a)]
    \item $W(t)$ is doubly stochastic.
    \item \label{cond:connecta} Every node in $\G(W(t))$ has a self-loop.
    \item \label{cond:bconnectivity}There exists an integer $B>0$ such that the  graph $\G_B(t)=([n],\EE_{B}(t))$ where
    \begin{align*}
       \EE_{B}(t)=\bigcup_{\tau=tB+1}^{(t+1)B}\EE(W(\tau))
    \end{align*}
    has a spanning rooted tree.
\end{enumerate}
\vspace{-0.3cm}
\QEDB
\end{assumption}

Finally, we assume the following standard condition on the step-size sequence $\{\eta(t)\}$.
\begin{assumption}\label{assump:AssumpOnStepsize}
The step-size $\eta$ satisfies $\eta(t)\geq 0$ for all $t\in\mathbb{Z}_{\geq0}$, $\sum_{t=0}^\infty \eta(t)=\infty$, and $\sum_{t=0}^\infty \eta^2(t)<\infty$.  \QEDB
\end{assumption}

\vspace{0.05cm}
The following proposition will play a key role in our analysis and proof of Theorem \ref{thm:MainTheorem}.

\vspace{0.05cm}
\begin{prop}\label{thm:AuxTheorem}
Suppose that Assumptions \ref{assump:AssumpOnFuncH}- \ref{assump:AssumpOnStepsize} hold. Then,  for all $i\in[n]$ and all initial conditions $\bx_i(0)\in \R$, the sequence $x_i$ generated by \eqref{eqn:AuxAlgorithm} satisfies
\begin{equation}
\lim_{t\to\infty}\bx_i(t)=z^\star,
\end{equation}
where $z^\star$ is a root of $H(z)$.  \QEDB
\end{prop}

%
%
The update rule \eqref{eqn:AuxAlgorithm} is inspired by distributed optimization algorithms that seek to find the roots of the gradient of a given cost function. However, in \eqref{eqn:AuxAlgorithm}, there is no global cost function to minimize, and, unlike the gradient of a convex function in optimization problems, in general, $H(z)$ is not non-decreasing. Moreover, note that  $h_i(z,t)$ is time-varying and might not be differentiable, which precludes the use of standard optimization techniques that rely on smoothness of the cost function and its derivatives \cite{popkov2005gradient,esteki2023distributed,simonetto2016class}. 
%
%


To prove Proposition \ref{thm:AuxTheorem}, we will also leverage the following deterministic variation of the Robbins-Siegmund Theorem \cite{ROBBINS1971233}:
\begin{prop}\label{thm:robsig}
Suppose that a non-negative sequence $\{{V}(t)\}$  satisfies  
\begin{align}\label{eqn:Lyapdecrease}
    {V}(t+1)\leq(1+{a}(t)){V}(t)-{b}(t)+{c}(t),
\end{align}
where ${a}(t),{b}(t),{c}(t)\geq 0$  for all $t$. If $\sum_{t=0}^{\infty}{a}(t)<\infty$ and $\sum_{t=0}^{\infty}{c}(t)<\infty$, then  $\lim_{t\to\infty}{V}(t)$ exists and $\sum_{t=0}^{\infty}{b}(t)<\infty$.   \QEDB
\end{prop}

Before proceeding, let's introduce a notational convention. For a vector $v = [v_1, \ldots, v_n]$, we denote the average of its elements as $\bar{v}$, defined as:
$$\bar{v} = \frac{1}{n} \sum_{j=1}^n v_j.$$
A key objective is for every region to determine the value $z^\star$. This necessitates that the differences between the state variables of any two regions vanish over time, i.e., $\lim_{t\to\infty} x_j(t) - x_k(t) = 0$ for all $j, k \in [n]$, or equivalently, $\lim_{t\to\infty} x_j(t) - \bar{x}(t) = 0$. To achieve this, we leverage a deterministic variation of a result presented in \cite[Lemma 9]{aghajan2022distributed}. 
\begin{lemma}[Lemma 9 in \cite{aghajan2022distributed}]\label{lem:Consensus}
Consider the update rule \eqref{eqn:AuxAlgorithm}. If $y_i(t)$ is bounded for all $i\in[n]$, and the sequence of matrices $\{W(t)\}$ satisfies Assumption \ref{assump:AssumpOnConnectivity}, then there exists a $\nu>0$ such that for all $i\in[n]$ and $t>0$
\begin{align*}
    |x_i(t)-\bar{x}(t)|\leq \nu\eta(t).
\end{align*}
for all $i\in[n]$ and $t>0$.  \QEDB
\end{lemma}

Note that Lemma \ref{lem:Consensus} not only ensures convergence but also quantifies the rate of convergence in terms of the sequence of step sizes $\eta$.

Lastly, we require the following lemma to establish the boundedness of $x(t)$.  The proof is presented in the Appendix.

\begin{lemma}\label{lem:BoundU}
Consider the discrete-time dynamics
\begin{align*}
        u(t+1)=u(t)-\eta(t)H(u(t))+w(t).
\end{align*}
where the sequences $\eta>0$,and $H$ are bounded. Suppose there exist $z^\star\in \R$ such that $(z-z^\star)H(z)\geq 0$ for all $z\in\R$. Moreover, assume the  sequence $\{w(t)\}$ is absolutely summable, i.e. $\sum_{t=1}^\infty |w(t)|\leq \Omega$ for some $\Omega>0$. Then $u(t)$ is bounded.  \QEDB
\end{lemma}

\vspace{0.1cm}
By leveraging the previous auxiliary results we can now establish Proposition \ref{thm:MainTheorem}.

{\it Proof of Proposition \ref{thm:AuxTheorem}:}
First, let us multiply both sides of \eqref{eqn:AuxAlgorithm} by $\frac{e^T}{n}$
\begin{align*}
    \frac{e^T}{n}{\bx}(t+1)=\frac{e^T}{n}W(t){\bx}(t)-\frac{e^T}{n}\eta(t){\by}(t).
\end{align*}
Since $W(t)$ is doubly stochastic,  we have
\begin{align*}
    \bar{{\bx}}(t+1)&=\bar{{\bx}}(t)-\eta(t){\bar{\by}}(t)\cr
    &\stackrel{(a)}{=}\bar{{\bx}}(t)-\frac{1}{n}\eta(t)\sum_{i=1}^n{\bh_i}(x_i(t))\cr
    &\quad+\frac{1}{n}\eta(t)\sum_{i=1}^n\left[{\bh_i}(x_i(t))-{\bh_i}(x_i(t),t)\right]\cr
    &\stackrel{(b)}{=}\bar{{\bx}}(t)-\frac{1}{n}\eta(t)\sum_{i=1}^n{\bh_i}(x_i(t))+\eta(t)A(t),
\end{align*}
where in $(a)$, $\bh_i(z)\triangleq\lim_{t\to\infty}h_i(z,t)$ for all $z\in \R$, and in $(b)$, $A(t)\triangleq \frac{1}{n}\sum_{i=1}^n\left[{\bh_i}(x_i(t))-{\bh_i}(x_i(t),t)\right]$.
Therefore, we have
\begin{align}\label{eqn:AveAuxAlg}
    \bar{{\bx}}(t+1)&=
    \bar{{\bx}}(t)-\eta(t)H(\bar{x}(t))\cr
    &\quad+\frac{1}{n}\eta(t)\sum_{i=1}^n\left[{\bh_i}(\bar{x}(t))-{\bh_i}(x_i(t))\right]+\eta(t)A(t)\cr
    &=
    \bar{{\bx}}(t)-\eta(t)H(\bar{x}(t))\cr
    &\quad+\frac{1}{n}\eta(t)\sum_{i=1}^n\Lambda_i(t)|\bar{x}(t)-x_i(t)|+\eta(t)A(t)\cr
    &=\bar{{\bx}}(t)-\eta(t)H(\bar{x}(t))+\eta(t)[A(t)+B(t)],
\end{align}
where $B(t)\triangleq \frac{1}{n}\sum_{i=1}^n\Lambda_i(t)|\bar{x}(t)-x_i(t)|$. Assumptions \ref{assump:AssumpOnFuncH}-2), \ref{assump:AssumpOnConvergenceFunction}, and \ref{assump:AssumpOnStepsize} imply that
\begin{align}\label{eqn:etaABsummable}
   \eta(t)[A(t)+B(t)]  \text{ is absolutely summable.}   
\end{align}
Therefore, by Lemma~\ref{lem:BoundU}, we have that $|\bar{x}(t)|$, and hence $|\bar{x}(t)-z^\star|$, are bounded, i.e., there exists $X>0$ such that for all $t$, we have $|\bar{x}(t)-z^\star|<X$. 

Let $V(t)\triangleq(x(t)-z^\star)^2$, where $z^\star$ is any point that satisfies Assumption \ref{assump:AssumpOnFuncH}-3). From \eqref{eqn:AveAuxAlg}, we can derive
\begin{align*}
    V(t+1)&=(\bar{{\bx}}(t+1)-z^\star)^2\cr
    &=\Big{(}\bar{{\bx}}(t)-z^\star-\eta(t)H(\bar{x}(t))+\eta(t)[A(t)+B(t)]\Big{)}^2\cr
    &=(\bar{{\bx}}(t)-z^\star)^2
    -2\eta(t)(\bar{{\bx}}(t)-z^\star)\bbh(\bar{x}(t))\cr
    &\quad+\eta^2(t)\bbh^2(\bar{x}(t))+2(\bar{{\bx}}(t)-z^\star)\eta(t)[A(t)+B(t)]\cr
    &\quad-2\eta^2(t)\bbh(\bar{x}(t))[A(t)+B(t)]\cr
    &\leq V(t)-b(t)+c(t),
\end{align*}
where
\begin{align*}
    b(t)&=2\eta(t)(\bar{{\bx}}(t)-z^\star)\bbh(\bar{x}(t)),\\
    c(t)&=\eta^2(t)\bbh^2(\bar{x}(t))+2|\bar{{\bx}}(t)-z^\star|\eta(t)|A(t)+B(t)|\cr
    &\quad+2|\bbh(\bar{x}(t)|\eta^2(t)|A(t)+B(t)|
\end{align*}
Note that Assumption \ref{assump:AssumpOnFuncH}-3) implies $b(t)\geq 0$. Assumption \ref{assump:AssumpOnFuncH}-1), \ref{assump:AssumpOnStepsize}, and \ref{eqn:etaABsummable} imply $\eta^2(t)\bbh^2(\bar{x}(t))+2|\bbh(\bar{x}(t))|\eta^2(t)|A(t)+B(t)|$ is  summable. Moreover, the fact that $|\bar{x}(t)-z^\star|$ is bounded and \ref{eqn:etaABsummable}  imply that $2|\bar{{\bx}}(t)-z^\star|\eta(t)|A(t)+B(t)|$ is  summable. Therefore, Proposition \ref{thm:robsig} implies 
\[\lim_{t\to\infty}(\bar{x}(t)-z^\star)^2 \text{ exists, and hence}, \lim_{t\to\infty}|\bar{x}(t)-z^\star| \text{ exists}\] 
and, moreover,
\begin{align}\label{eqn:btsummable}
    \sum_{t=0}^\infty \eta(t)(\bar{{\bx}}(t)-z^\star)\bbh(\bar{x}(t))<\infty.
\end{align}
Next, we show that $\bar{x}(t)$ converges to a root of $H(z)$. If $\lim_{t\to\infty}|\bar{x}(t)-z^\star|=0$, there is nothing to prove. Hence, assume that  $\lim_{t\to\infty}|\bar{x}(t)-z^\star|>0$. From Assumption \ref{assump:AssumpOnFuncH}-3), we have $(\bar{{\bx}}(t)-z^\star)\bbh(\bar{x}(t))\geq 0$, which implies 
\begin{align*}
    (\bar{{\bx}}(t)-z^\star)\bbh(\bar{x}(t))=|\bar{{\bx}}(t)-z^\star||\bbh(\bar{x}(t))|.
\end{align*}
Therefore, \eqref{eqn:btsummable} and the fact that $\sum_{t=0}^\infty \eta(t)=\infty$ imply $\liminf_{t\to\infty}|\bbh(\bar{x}(t))|=0$, and hence 
$\liminf_{t\to\infty}\bbh(\bar{x}(t))=0$. Therefore, there exists a sub-sequence $\{H(\bar{x}(t_k))\}_{k=0}^\infty$ such that $\lim_{k\to\infty} H(\bar{x}(t_k))=0$. Since $H(\cdot)$ is a continuous function according to Assumption \ref{assump:AssumpOnFuncH}-2, we have 
\[ H\left(\lim_{k\to\infty}\bar{x}(t_k)\right)=\lim_{k\to\infty} H(\bar{x}(t_k))=0.\]
Hence $\tilde{z}=\lim_{k\to\infty}\bar{x}(t_k)$ is a root of $H(z)$. Since $z^\star$ can be any point that satisfies Assumption \ref{assump:AssumpOnFuncH}-3), we set $z^\star=\tilde{z}$. Hence $\lim_{t\to\infty}|\bar{x}(t)-z^\star|=0$. 

Finally, we have 
\begin{align*}
    \lim_{t\to\infty}|x_i(t)-z^\star|\leq\lim_{t\to\infty}|\bar{x}(t)-z^\star|+\lim_{t\to\infty}|x_i(t)-\bar{x}(t)|=0,
\end{align*}
which follows by Lemma \ref{lem:Consensus}. This completes the proof. \hfill $\blacksquare$
\subsection{ Proof of Theorem \ref{thm:MainTheorem}} 
Now, we are ready to prove the main result of the paper, Theorem \ref{thm:MainTheorem}. First, we need to show in Algorithm \ref{alg:DisLoadShed} that $x_j(t)$ converges to $\hat{z}$, for all $j\in[n]$, where $\hat{f}(\hat{z})=\P$. Thus, we proceed to verify that all the assumptions of Proposition \ref{thm:AuxTheorem} are satisfied.

To construct a doubly stochastic sequence $\{W(t)\}$, we use the Metropolis-Hastings algorithm. In particular, since the communication network satisfies Assumption \ref{assump:Network}, the matrix sequence $\{W(t)\}$ satisfies Assumption \ref{assump:AssumpOnConnectivity}-(c) \cite{tahbaz2006one}. Similarly, by design, we construct a step-size sequence that satisfies Assumption \ref{assump:AssumpOnStepsize}.  Therefore, it suffices to verify whether $\hat{f}_j(z)-p_j(t)$s and $\hat{f}(z)-\P$ satisfy Assumption \ref{assump:AssumpOnFuncH}. Since the number of loads and their criticality values are finite, $\hat{f}_j(z)$s and $\hat{f}(z)$ are bounded. Hence, due to the boundness of $p_j(t)$,  Assumption \ref{assump:AssumpOnFuncH}-1) is satisfied. Additionally, since $\bw(z)$ is Lipschitz,  $\hat{f}_j(z)$ and $\hat{f}(z)$ are Lipschitz and satisfy Assumption \ref{assump:AssumpOnFuncH}-2). Finally, we have that
\begin{equation}
 \lim_{z\to\infty}\hat{f}(z)-\P=\sum_{\ell\in\L}\ell- \P\geq 0,
\end{equation}
due to Assumption \ref{assump:sumellgeP} and $\hat{f}(0)-\P\leq 0$. Since $\hat{f}(z)$ is a continuous non-decreasing function we also have that $\hat{f}(z)-\P$ satisfies Assumption \ref{assump:AssumpOnFuncH}-3. Therefore, $x_j(t)$ converges to $\hat{z}$ for all $j\in[n]$.

To complete the proof, we show that $z_j(t)$ converges to $z^*$ for all $j\in[n]$ in finite time.
To do so, first set 
\begin{align}\label{eqn:defEps}
    \epsilon\triangleq\frac{1}{2}\min\left\{\hat{z}-\max_{\ell\in\L:C(\ell)\leq \hat{z}}C(\ell),z^*-\hat{z}\right\}.
\end{align}
Note that since $f(z^*)>\P$, from Lemma \ref{lem:ApproxZstar}, we have 
$z^*=\min_{\ell\in\L:C(\ell)> \hat{z}}C(\ell)$.
Hence, $z^*>\hat{z}$. Moreover, due to \eqref{eqn:defEps}, we have 
\begin{align}\label{eqn:zStarzhatEps}
    z^*>\hat{z}+\epsilon.
\end{align}
On the other hand, $\hat{z}>z'\triangleq\max_{\ell\in\L:C(\ell)\leq \hat{z}}C(\ell)$, since otherwise, $f(z')\geq \P$, which contradicts the fact that $f(z)<\P$ for all $z<z^*$.
This fact, together with  $z^*>\hat{z}$ imply that $\epsilon>0$. Additionally, due to \eqref{eqn:defEps}, we have 
\begin{align}\label{eqn:zhatMaxEps}
    \hat{z}>\max_{\ell\in\L:C(\ell)\leq \hat{z}}C(\ell)+\epsilon.
\end{align}
Using \eqref{eqn:zStarzhatEps} and \eqref{eqn:zhatMaxEps}, we obtain
\begin{align}\label{eqn:ThereIsNo}
    \text{there is no $\ell\in\L$ such that $\hat{z}-\epsilon\leq C(\ell)\leq \hat{z}+\epsilon$}.
\end{align}
Since $x_j(t)\rightarrow\hat{z}$ for all $j\in[n]$, there exist $T_\epsilon$ such that $|x_j(t)-\hat{z}|<\epsilon$ for all $j\in [n]$ and $t>T_\epsilon$. Therefore, for $t>T_\epsilon$, we have
\begin{align*}
    \zeta_j(t)&=\min_{\ell\in \L_j :C(\ell)\geq x_j(t)} C(\ell)\cr
    &=\min_{\ell\in \L_j :C(\ell)\geq \hat{z}+x_j(t)-\hat{z}} C(\ell)\cr
    &\stackrel{(a)}{=}\min_{\ell\in \L_j :C(\ell)\geq \hat{z}} C(\ell)\triangleq\hat{\zeta}_j,
\end{align*}
where $(a)$ follows from \eqref{eqn:ThereIsNo}.
Therefore, for $t>T_\epsilon$, $\zeta_j(t)$ for all $j$ is constant, which implies convergence in finite time. \hfill $\blacksquare$
%
\section{EXPERIMENTAL RESULTS}
\label{sec:simulation}
In this section, we first present the model used to implement our algorithm in the Quebec 29-bus system. We then explore a case study focused on a specific communication network and cumulative criticality functions, applying our algorithm and discussing the results. Lastly, we showcase the flexibility of our approach by extending its application to the continuous variation of the load-shedding problem discussed in Section V.B of \cite{fitri2022priority}.

\subsection{Use Case}
\label{sec3}

\begin{figure}[t!]
\centering
\includegraphics[width=0.45\textwidth]{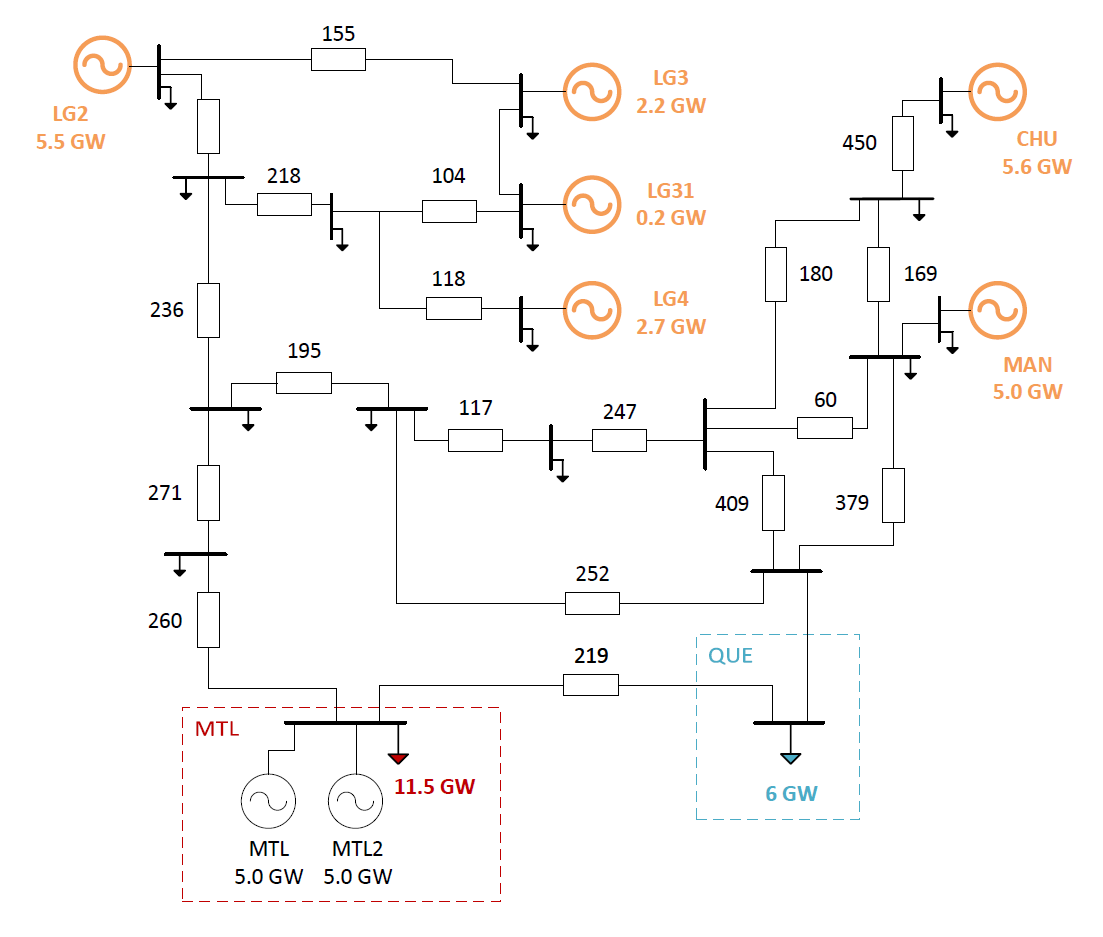}
\caption{\small{Quebec 29-bus system diagram.}}
\centering
\label{fig:Quebec29}
 \vspace{-0.5cm}
\end{figure}

The proposed algorithm was tested on the Quebec 29-bus system, a simplified model of Quebec's transmission power grid. This model, available on the MathWorks website, features two generation areas in the North and two load centers (MTL and QUE) in the South \cite{Crivellaro2019BeyondLS}. Power is transported via 735 kV lines. The system, dominated by hydropower, has 8 generating units with a total capacity of 26.2 GW.
Figure \ref{fig:Quebec29} depicts the system diagram. Generators disconnected for testing the adaptive UFLS algorithm are highlighted in orange, with annotations indicating their names and ratings. Aggregated loads for MTL and QUE are shown in red and blue, respectively, with their capacities specified. Black arrows represent non-sheddable loads, and line lengths are given in kilometers.
Quebec's frequency requirements, due to its hydropower-dominant system, are more flexible than those in the US, as detailed in NERC's PRC-024-3 standard \cite{NorthAmericanElectricReliabilityCorporation2020}.

\subsection{Numerical Results}
\begin{figure}[t!]
\centering
\begin{tikzpicture}[
  every node/.style={draw, circle},  
  minimum size=1cm,
  thick]                

\node (1) at (0,0) {MTL1};
\node (2) at (2,0) {MTL2};
\node (3) at (4,0) {QUE1};
\node (4) at (6,0) {QUE2};

\draw[<->] (1) -- (2);  
\draw[<->] (2) -- (3); 
\draw[<->] (3) -- (4);
\end{tikzpicture}
\caption{\small{Communication Network: Illustrating both In-Area and Out-of-Area connections between four regions (MTL1, MTL2, QUE1, QUE2).}}
\label{fig:line_graph}
\end{figure}
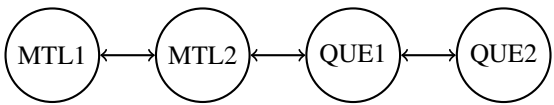

\begin{figure}[t!]
\centering
\includegraphics[width=0.49\textwidth,height=0.65\linewidth]{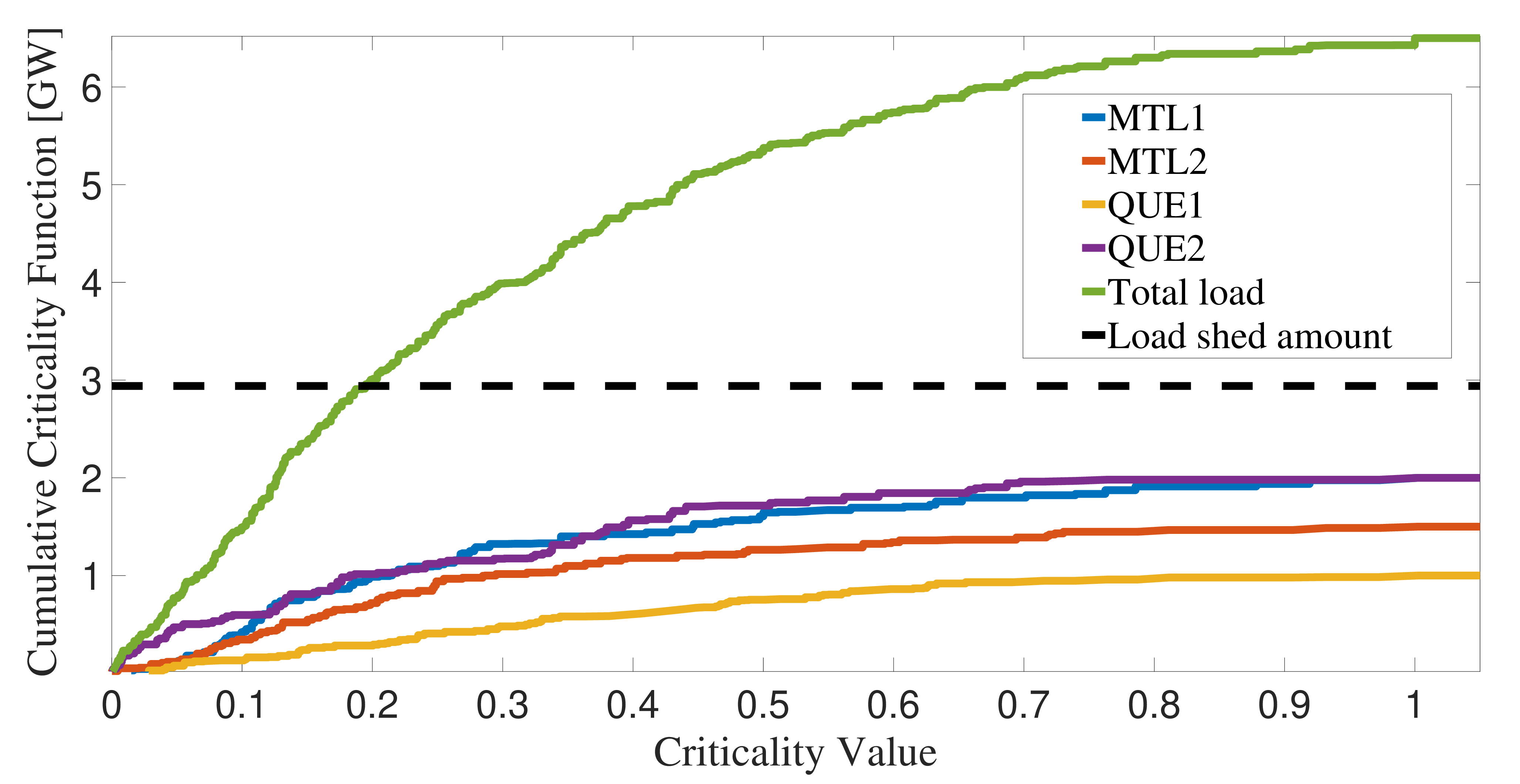}
\caption{\small{Cumulative criticality function of four regions MTL1, MTL2, QUE1, and QUE2, and total amount of loads. }}
\centering
\label{fig:CCFSim}
 \vspace{-0.5cm}
\end{figure}

In our simulation model, the main areas of Montreal and Québec are each subdivided into two regions MTL1, MTL2, QUE1, and QUE2, respectively. The resulting four regions communicate via a network that allows information exchange every 10 milliseconds whenever a communication link exists (or is available) between them.
The communication network is illustrated in Figure \ref{fig:line_graph}. 
For instance, MTL1 and QUE2, despite being part of the overall system, cannot directly communicate with each other.
Each region contains 100 loads, and the criticality value of each load is known only within its respective region. Figure \ref{fig:CCFSim} illustrates the corresponding CCFs for each region.

In our simulation, a 2.7 GW power loss occurs at $t=4$s due to the LG4 generator tripping. To stabilize the grid, a 2.94 GW load shed is needed. However, the initial instantaneous precise amount of load shedding is assumed to be unknown. In particular, the load shedding in each of the four regions is given by
  $p_i(t) = 2.94 \times \frac{1}{4} + \frac{e_i(t)}{t}$ GW, 
where $i$ represents the region, and $e_i(t)$ is a noise term with a uniform distribution between $-1$ and $1$, which can capture uncertainty from local estimation errors based on the ROCOF.

Load shedding is delayed until the frequency drops below 59.5 Hz at $t=5.05$s.
During this critical interval (4s to 5.05s), \textbf{Algorithm 1} is used to determine which loads to shed based on their criticality. 
The algorithm employs a Metropolis-Hastings matrix and step-size $\eta(t) = \frac{1}{t+1}$.
In our simulation, we also incorporate a dynamic min-consensus protocol with local self-tuning \cite{deplano2022dynamicss} running concurrently with Algorithm \ref{alg:DisLoadShed}.  Specifically, we simultaneously execute the following dynamics:
\begin{align*}
    z_j(t+1)&\leftarrow\min_{k\in \N_j\cup\{j\}}\left\{z_k(t)+\alpha_j(t),\zeta_j(t+1)\right\}\\
    \alpha_j(t+1)&\leftarrow\begin{cases}
        \frac{1}{2},& \text{if }z_j(t+1)>z_j(t)\\
        \frac{c}{2},& \text{otherwise}
    \end{cases}.
\end{align*}

\begin{figure}[h!]
\includegraphics[width=0.53\textwidth,height=0.74\linewidth]{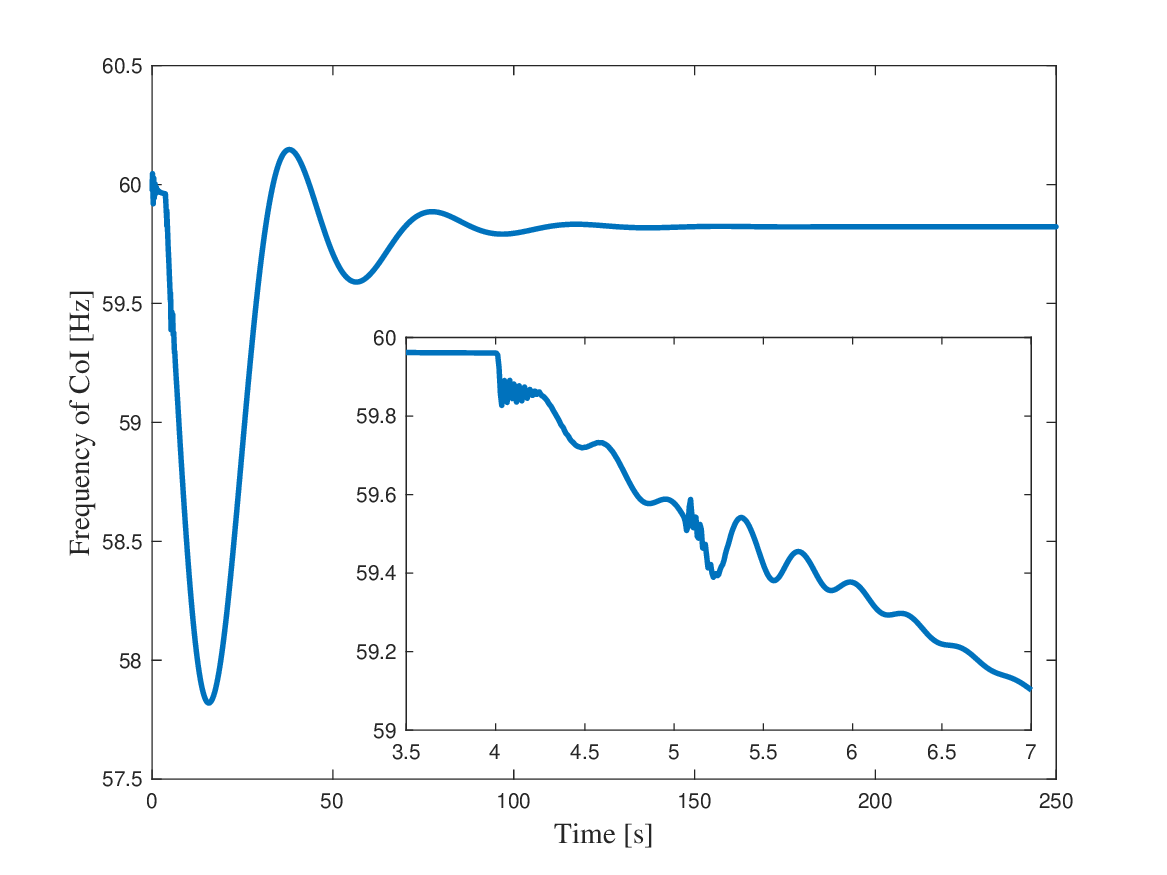}
\centering
\caption{\small{Frequency at the system's center of inertia (CoI) after disconnecting the generator LG4 and load shedding. The smaller box provides a magnified view of the plot.}
}
\label{fig:CoILG4}
 \vspace{-0.5cm}
\end{figure}

\begin{figure}[h!]
\includegraphics[width=0.53\textwidth,height=0.74\linewidth]{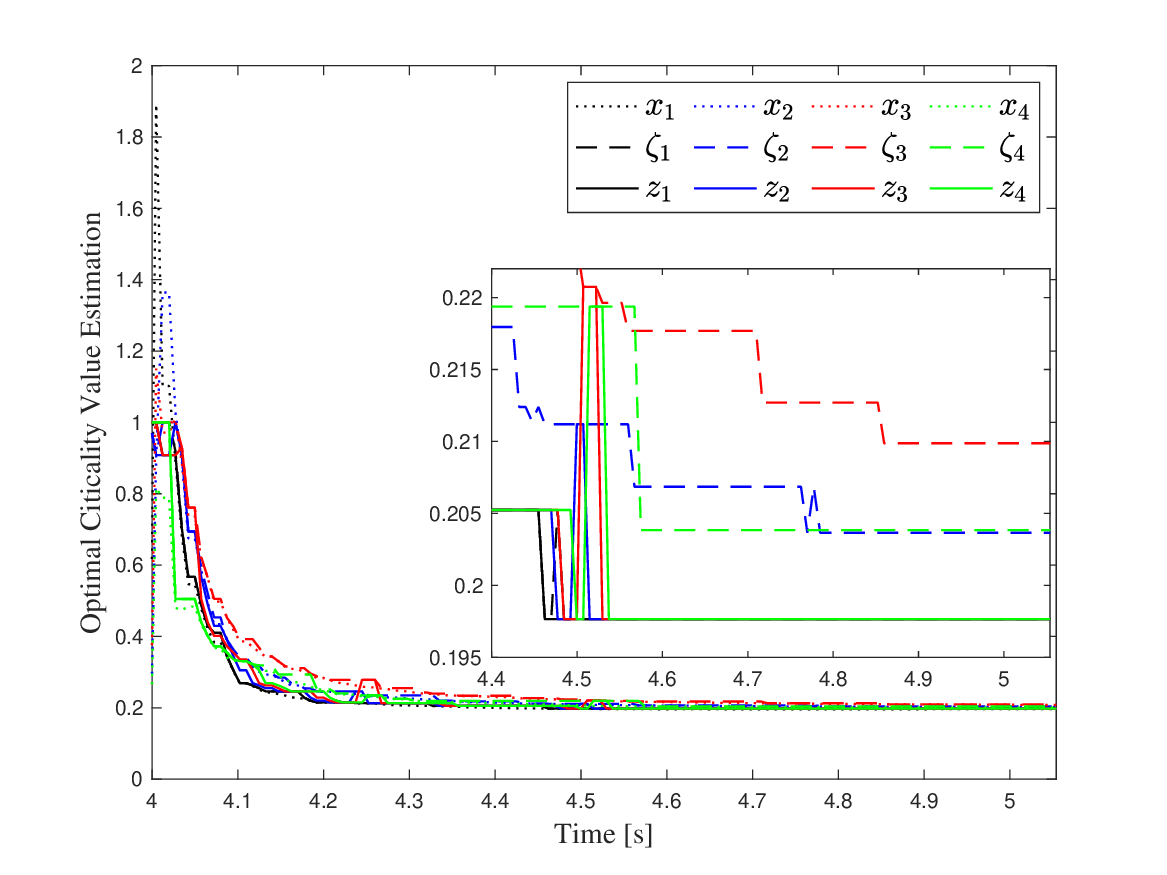}
\centering
\caption{\small{Estimation of optimal criticality value in four regions with distributed (line) network. The smaller box provides a magnified view of the plot.
}}
\label{fig:OptimalEst}
 \vspace{-0.5cm}
\end{figure}

Figure \ref{fig:CoILG4} shows the algorithm's success in stabilizing the system: the center of inertia frequency drops sharply after the trip but quickly recovers once load shedding is initiated.
Figure \ref{fig:OptimalEst} shows the agent's estimation of the optimal criticality value given by 0.1947. After 5 seconds, the algorithm converges to 0.1977.  The optimal load-shedding amount for each agent is given by 0.9585 GW for MTL1, 0.6966 GW for MTL2, 0.2834 GW for QUE1, and 1.0147 GW for QUE2.

\subsection{Application to the Continuous Variation Case}
\begin{figure}[t!]
\includegraphics[width=0.53\textwidth,height=0.74\linewidth]{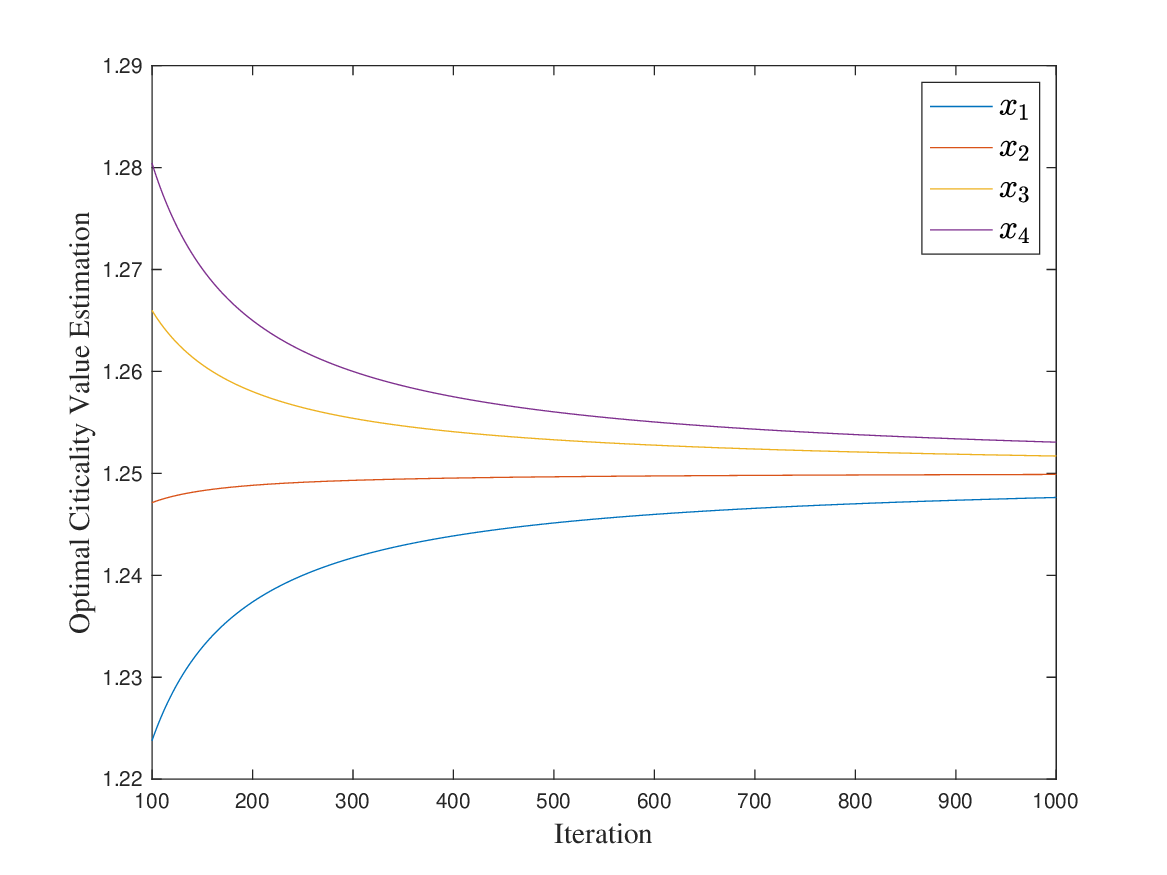}
\centering
\caption{\small{Estimation of optimal criticality for the example in \cite{fitri2022priority}.}
}
\label{fig:OptimalEstContinuous}
\vspace{-0.3cm}
\end{figure}

To further showcase the versatility of the proposed approach, we apply \textbf{Algorithm 1} to the continuous load-shedding problem investigated in \cite{fitri2022priority}. In particular, for the purpose of comparison, we consider the same illustrative example presented in Section V.B of \cite{fitri2022priority}. We consider a system with four regions, each having the capacity to shed a continuous load amount up to a maximum of 1.2 GW. Since in this case, all loads are assumed to be of the same nature (i.e., $C_n$ is immaterial), the criticality values are completely characterized by the regions, and they are defined as follows: $C_r(1) = 1$, $C_r(2) = C_r(3) = 2$, and $C_r(4) = 3$. The objective is to shed a total of 1.8 GW of load, prioritizing regions with lower criticality values. We apply our algorithm to this scenario and show the resulting trajectories converging to the estimated optimal criticality value in Figure \ref{fig:OptimalEstContinuous}. The results demonstrate that after 1000 iterations, the criticality estimation converges to the values $[1.2477, 1.2507, 1.2532, 1.2550]$. Using the solution approach described in Equation \eqref{eqn:ContinuousSol}, this leads to the following load-shedding distribution across the four regions: $[1.2,0.3009,0.3038,0]$ GW. This successful application of our algorithm to the continuous load-shedding problem highlights its adaptability and effectiveness in handling scenarios where load shedding can be performed in continuous spaces.

\section{CONCLUSIONS}\label{Conclusion}
This paper presents a novel approach to distributed resource allocation problems on discrete action spaces, with applications to priority-based distributed load shedding in power systems. The proposed approach addresses limitations of traditional resource allocation algorithms by considering discrete load values, heterogeneous criticality within regions, and time-varying communication networks. By leveraging the concept of \emph{cumulative criticality functions} (CCFs), the optimal resource allocation problem can be cast as finding the root of a time-varying function, enabling a distributed solution. Through rigorous analysis, we established the convergence of our algorithm under dynamic network conditions. Numerical simulations on the Quebec 29-bus system validated the effectiveness of our approach in achieving adaptive under-frequency load shedding while adhering to predefined priority orders and adapting to network changes. While the approach was presented in the context of load shedding in energy systems, it can also be applied to other optimal resource allocation problems defined on discrete spaces, where traditional gradient-based methods are not applicable.

\appendix
{\it Proof of Lemma \ref{obs:zstar}:} Let $\tilde{z}\in \R$ such that $\tilde{z}\not\in\{C(\ell)|\ell\in\L\}$, and let 
$$z'=\min_{z\in \{C(\ell)|\ell\in\L\}:z<\tilde{z}}z.$$
Then, from the definition of CCF, we have $f(\tilde{z})=f(z')$ where $z'<\tilde{z}$, which implies $\tilde{z}$ cannot be $z^*$ defined in \eqref{eqn:Sol2}. \hfill $\blacksquare$

\vspace{0.1cm}
{\it Proof of Lemma  \ref{obs:fequalfhat}:} 
Since $c\leq C(\ell')-C(\ell)$, for all $\ell'$ such that $C(\ell')>C(\ell)$,
we have ${\bw_c(C(\ell)-C(\ell'))=0}$ for all of those $\ell'$. Hence, $\bw_c(z)=1$ for $z\geq 0$ implies $\hat{f}(C(\ell))=\sum_{\tilde{\ell}\in\L:C(\tilde{\ell})\leq C(\ell)}\tilde{\ell}$, which is equal to $f(C(\ell))$.  \hfill $\blacksquare$

{\it Proof of Lemma \ref{lem:ApproxZstar}:} 
First, from Lemma \ref{obs:zstar} and \ref{obs:fequalfhat}, we have $f(z^*)=\hat{f}(z^*)$. Moreover, there is no $\tilde{\ell}\in\L$, such that $\min\{z^*,\hat{z}\}<C(\tilde{\ell})<\max\{z^*,\hat{z}\}$. 
This is because 
if $\hat{z}<C(\tilde{\ell})<z^*$, then $f(C(\tilde{\ell}))\geq \P$ which contradicts the fact that $f(z)<\P$ for all $z<z^*$. On the other hand, if $z^*<C(\tilde{\ell})<\hat{z}$, then due to $z^*<C(\tilde{\ell})$, we have $f(C(\tilde{\ell}))>f(z^*)\geq\P$, which contradicts the fact that $f(C(\tilde{\ell}))=\hat{f}(C(\tilde{\ell}))\leq\P$ due to $C(\tilde{\ell})<\hat{z}$.
Finally, the facts that $f(z^*)>\P$ implies $z^*> \hat{z}$ and $f(z^*)=\P$ implies $z^*\leq \hat{z}$ complete the proof.
 \hfill $\blacksquare$

\vspace{0.1cm}
 {\it Proof of Lemma \ref{lem:BoundU}:} First note that since $(z-z^\star)H(z)\geq 0$ for all $z\in\R$, we have $u(t)\geq z^\star$ if $H(u(t))\geq 0$, and $u(t)\leq z^\star$ if $H(u(t))\leq 0$. Therefore, we only need to prove $u(t)$ has a lower (upper) bound, if $H(u(t))\leq 0$ ($H(u(t))\geq 0$).
Let $H(u(t))\leq 0$ and  
\[T(t)\triangleq\min_{H(u(k))\leq 0 \text{ for } \tau\leq k\leq t}\tau.\]
We have
\begin{align*}
    u(t)&=u(T(t))-\sum_{k=T(t)}^{t-1}\eta(k)H(u(k))+\sum_{k=T(t)}^{t-1} w(k)\cr
    &\geq u(T(t))-\sum_{k=0}^{\infty} |w(k)|\geq u(T(t))-\Omega.
\end{align*}
If $T(t)=0$, then $u(t)$ has a lower bound. If not, we have
\begin{align*}
    u(t)&\geq
    u(T(t))-\Omega\cr
    &=u(T(t)-1)-\eta(T(t)-1)H(x(T(t)-1))\cr
    &\quad+w(x(T(t)-1))-\Omega\cr
    &\geq z^\star- \eta(T(t)-1)H(x(T(t)-1))+w(x(T(t)-1))-\Omega,
\end{align*}
where follows from the fact that $u(T(t)-1)>0$. Therefore, the fact that $\eta(t)$, $w(t)$, and $H(z)$ are bounded implies $u(t)$ has lower bound.
 Similarly, we can prove $u(t)$ has an upper bound if $H(u(t))\geq 0$. \hfill $\blacksquare$

\vspace{-0.4cm}
\bibliographystyle{IEEEtran}
\bibliography{sample.bib}

\end{document}